\documentclass[12pt]{amsart}

\usepackage[margin=3cm]{geometry}

\usepackage{amsmath,amssymb,amsbsy,amsfonts,amsthm,latexsym,
                        amsopn,amstext,amsxtra,euscript,amscd,mathrsfs,color,bm, cite}
                        \usepackage{ulem}

\usepackage{float}
\usepackage[english]{babel}
\usepackage{mathtools}
\usepackage{todonotes}
\usepackage{url}
\usepackage[colorlinks,linkcolor=blue,anchorcolor=blue,citecolor=blue,backref=page]{hyperref}

\usepackage[norefs,nocites]{refcheck}

\newtheorem{theorem}{Theorem}
\newtheorem{lemma}[theorem]{Lemma}

\newtheorem{cor}[theorem]{Corollary}

\numberwithin{equation}{section}
\numberwithin{theorem}{section}
\numberwithin{table}{section}

\numberwithin{figure}{section}

\begin{document}


\newfont{\teneufm}{eufm10}
\newfont{\seveneufm}{eufm7}
\newfont{\fiveeufm}{eufm5}
%
%
\newfam\eufmfam
                \textfont\eufmfam=\teneufm \scriptfont\eufmfam=\seveneufm
                \scriptscriptfont\eufmfam=\fiveeufm
%
%
\def\frak#1{{\fam\eufmfam\relax#1}}
%


\def\bbbr{{\rm I\!R}} 
\def\bbbm{{\rm I\!M}}
\def\bbbn{{\rm I\!N}} 
\def\bbbf{{\rm I\!F}}
\def\bbbh{{\rm I\!H}}
\def\bbbk{{\rm I\!K}}
\def\bbbp{{\rm I\!P}}
\def\bbbone{{\mathchoice {\rm 1\mskip-4mu l} {\rm 1\mskip-4mu l}
{\rm 1\mskip-4.5mu l} {\rm 1\mskip-5mu l}}}
\def\bbbc{{\mathchoice {\setbox0=\hbox{$\displaystyle\rm C$}\hbox{\hbox
to0pt{\kern0.4\wd0\vrule height0.9\ht0\hss}\box0}}
{\setbox0=\hbox{$\textstyle\rm C$}\hbox{\hbox
to0pt{\kern0.4\wd0\vrule height0.9\ht0\hss}\box0}}
{\setbox0=\hbox{$\scriptstyle\rm C$}\hbox{\hbox
to0pt{\kern0.4\wd0\vrule height0.9\ht0\hss}\box0}}
{\setbox0=\hbox{$\scriptscriptstyle\rm C$}\hbox{\hbox
to0pt{\kern0.4\wd0\vrule height0.9\ht0\hss}\box0}}}}
\def\bbbq{{\mathchoice {\setbox0=\hbox{$\displaystyle\rm
Q$}\hbox{\raise 0.15\ht0\hbox to0pt{\kern0.4\wd0\vrule
height0.8\ht0\hss}\box0}} {\setbox0=\hbox{$\textstyle\rm
Q$}\hbox{\raise 0.15\ht0\hbox to0pt{\kern0.4\wd0\vrule
height0.8\ht0\hss}\box0}} {\setbox0=\hbox{$\scriptstyle\rm
Q$}\hbox{\raise 0.15\ht0\hbox to0pt{\kern0.4\wd0\vrule
height0.7\ht0\hss}\box0}} {\setbox0=\hbox{$\scriptscriptstyle\rm
Q$}\hbox{\raise 0.15\ht0\hbox to0pt{\kern0.4\wd0\vrule
height0.7\ht0\hss}\box0}}}}
\def\bbbt{{\mathchoice {\setbox0=\hbox{$\displaystyle\rm
T$}\hbox{\hbox to0pt{\kern0.3\wd0\vrule height0.9\ht0\hss}\box0}}
{\setbox0=\hbox{$\textstyle\rm T$}\hbox{\hbox
to0pt{\kern0.3\wd0\vrule height0.9\ht0\hss}\box0}}
{\setbox0=\hbox{$\scriptstyle\rm T$}\hbox{\hbox
to0pt{\kern0.3\wd0\vrule height0.9\ht0\hss}\box0}}
{\setbox0=\hbox{$\scriptscriptstyle\rm T$}\hbox{\hbox
to0pt{\kern0.3\wd0\vrule height0.9\ht0\hss}\box0}}}}
\def\bbbs{{\mathchoice
{\setbox0=\hbox{$\displaystyle     \rm S$}\hbox{\raise0.5\ht0\hbox
to0pt{\kern0.35\wd0\vrule height0.45\ht0\hss}\hbox
to0pt{\kern0.55\wd0\vrule height0.5\ht0\hss}\box0}}
{\setbox0=\hbox{$\textstyle        \rm S$}\hbox{\raise0.5\ht0\hbox
to0pt{\kern0.35\wd0\vrule height0.45\ht0\hss}\hbox
to0pt{\kern0.55\wd0\vrule height0.5\ht0\hss}\box0}}
{\setbox0=\hbox{$\scriptstyle      \rm S$}\hbox{\raise0.5\ht0\hbox
to0pt{\kern0.35\wd0\vrule height0.45\ht0\hss}\raise0.05\ht0\hbox
to0pt{\kern0.5\wd0\vrule height0.45\ht0\hss}\box0}}
{\setbox0=\hbox{$\scriptscriptstyle\rm S$}\hbox{\raise0.5\ht0\hbox
to0pt{\kern0.4\wd0\vrule height0.45\ht0\hss}\raise0.05\ht0\hbox
to0pt{\kern0.55\wd0\vrule height0.45\ht0\hss}\box0}}}}
\def\bbbz{{\mathchoice {\hbox{$\sf\textstyle Z\kern-0.4em Z$}}
{\hbox{$\sf\textstyle Z\kern-0.4em Z$}} {\hbox{$\sf\scriptstyle
Z\kern-0.3em Z$}} {\hbox{$\sf\scriptscriptstyle Z\kern-0.2em
Z$}}}}
\def\ts{\thinspace}

\def\squareforqed{\hbox{\rlap{$\sqcap$}$\sqcup$}}
\def\qed{\ifmmode\squareforqed\else{\unskip\nobreak\hfil
\penalty50\hskip1em\null\nobreak\hfil\squareforqed
\parfillskip=0pt\finalhyphendemerits=0\endgraf}\fi}

\def\cA{{\mathcal A}}
\def\cB{{\mathcal B}}
\def\cC{{\mathcal C}}
\def\cD{{\mathcal D}}
\def\cE{{\mathcal E}}
\def\cF{{\mathcal F}}
\def\cG{{\mathcal G}}
\def\cH{{\mathcal H}}
\def\cI{{\mathcal I}}
\def\cJ{{\mathcal J}}
\def\cK{{\mathcal K}}
\def\cL{{\mathcal L}}
\def\cM{{\mathcal M}}
\def\cN{{\mathcal N}}
\def\cO{{\mathcal O}}
\def\cP{{\mathcal P}}
\def\cQ{{\mathcal Q}}
\def\cR{{\mathcal R}}
\def\cS{{\mathcal S}}
\def\cT{{\mathcal T}}
\def\cU{{\mathcal U}}
\def\cV{{\mathcal V}}
\def\cW{{\mathcal W}}
\def\cX{{\mathcal X}}
\def\cY{{\mathcal Y}}
\def\cZ{{\mathcal Z}}

\def\HH{{\mathsf H}}

\def\le{\leqslant}
\def\leq{\leqslant}
\def\ge{\geqslant}
\def\leq{\leqslant}

\def\sfB{\mathsf {B}}
\def\sfC{\mathsf {C}}
\def\sfS{\mathsf {S}}
\def\sfI{\mathsf {I}}
\def\sfT{\mathsf {T}}
\def\L{\mathsf {L}}
\def\FF{\mathsf {F}}

\def\sE {\mathscr{E}}
\def\sS {\mathscr{S}}

\def\nrp#1{\left\|#1\right\|_p}
\def\nrq#1{\left\|#1\right\|_m}
\def\nrqk#1{\left\|#1\right\|_{m_k}}
\def\Ln#1{\mbox{\rm {Ln}}\,#1}
\def\nd{\hspace{-1.2mm}}
\def\ord{{\mathrm{ord}}}

\def\va{{\mathbf{a}}}
\def\vb{{\mathbf{b}}}
\def\vc{{\mathbf{c}}}
\def\vn{{\mathbf{n}}}
\def\vr{{\mathbf{r}}}
\def\vu{{\mathbf{u}}}

\def \fB{\mathfrak{B}}
\def \fR{\mathfrak{R}}
\def \fU{\mathfrak{U}}
\def \fW{\mathfrak{W}}

\newcommand{\comm}[1]{\marginpar{%
\vskip-\baselineskip 
\raggedright\footnotesize
\itshape\hrule\smallskip#1\par\smallskip\hrule}}




\newcommand{\ignore}[1]{}

\def\vec#1{\mathbf{#1}}

\def\e{\mathbf{e}}
\def\em{\mathbf{e}_m}



\def\GL{\mathrm{GL}}

\hyphenation{re-pub-lished}

\def\rank{{\mathrm{rk}\,}}
\def\dd{{\mathrm{dyndeg}\,}}
\def\xbar{\overline{x}}
\def\sbar{\overline{s}}

\def\A{\mathbb{A}}
\def\B{\mathbf{B}}
\def \C{\mathbb{C}}
\def \F{\mathbb{F}}
\def \K{\mathbb{K}}

\def \Z{\mathbb{Z}}
\def \P{\mathbb{P}}
\def \R{\mathbb{R}}
\def \Q{\mathbb{Q}}
\def \N{\mathbb{N}}
\def \Z{\mathbb{Z}}

\def \balpha{\bm{\alpha}}
\def \bbeta{\bm{\beta}}
\def \bgamma{\bm{\gamma}}
\def \blambda{\bm{\lambda}}
\def \bchi{\bm{\chi}}
\def \bphi{\bm{\varphi}}
\def \bpsi{\bm{\psi}}
\def \bomega{\bm{\omega}}
\def \btheta{\bm{\vartheta}}
\def \btheta{\bm{\vartheta}}

\def \L{\mathbf{L}}

\def \nd{{\, | \hspace{-1.5 mm}/\,}}

\def\mand{\qquad\mbox{and}\qquad}

\def\Zn{\Z_n}

\def\Fp{\F_p}
\def\Fq{\F_q}
\def \fp{\Fp^*}
\def\\{\cr}
\def\({\left(}
\def\){\right)}
\def\fl#1{\left\lfloor#1\right\rfloor}
\def\rf#1{\left\lceil#1\right\rceil}
\def\vh{\mathbf{h}}

\def\e{{\mathbf{\,e}}}
\def\ep{{\mathbf{\,e}}_p}

\def\eps{\varepsilon}

\newcommand{\lcm}{\operatorname{lcm}}
\newcommand{\ind}{\operatorname{ind}}
\newcommand{\Res}{\operatorname{Res}}

\def\cFH{\mathcal{F}_{k}(H)}
\def\cGH{\mathcal{G}_{g}(H)}

\def\sssum{\mathop{\sum\!\sum\!\sum}}
\def\ssum{\mathop{\sum\ldots \sum}}
\def\dsum{\mathop{\quad \sum \ \sum}}
\def\iint{\mathop{\int\ldots \int}}

\def\UkH{U_{k}(m,H,N)}
\def\WgH{W_{g}(m,H,N)}
\def\BH{\cB_{k}(H)}
\def\IH{\cI_g(H)}
\def \MbH{M_\vec{b}(m,H)}
\def \Tbm{T_\vec{b}(m)}
\def \Tbd{T_\vec{b}(m/d)}

\title[Congruences with factorials]{Additive congruences with factorials modulo a prime}

\author[M.~Z.~Garaev]{Moubariz~Z.~Garaev}  \author[J. C. Pardo] {Julio C. Pardo}
\address{Centro  de Ciencias Matem{\'a}ticas,  Universidad Nacional Aut\'onoma de
M{\'e}\-xico, C.P. 58089, Morelia, Michoac{\'a}n, M{\'e}xico}
\email{garaev@matmor.unam.mx}

\email{jcpardo@matmor.unam.mx}

\begin{abstract} 
Let $p$ be a large prime number. We prove that any integer $\lambda$ modulo $p$ can be represented in the form
$$
m!n! +\sum_{i=1}^{47}n_i!\equiv \lambda \pmod p,
$$
with $\max\{m,n,n_1,\ldots,n_{47}\}\ll p^{1300/1301}.$ This improves the exponent $1350/1351$ of Garaev, Luca and Shparlinski (2005). Furthermore, we prove that any integer $\lambda$ can be represented in the form
$$
m_1!n_1! +m_2!n_2!+m_3!n_3! +m_4!n_4!+m_5!n_5! \equiv \lambda \pmod p
$$ 
with  $\max\{m_1,n_1,\ldots,m_5,n_5\}\le p^{97/113 +o(1)}.$ This improves the exponent $27/28$ of Garaev and Garcia (2007). The proofs of these two results are based on the recent work of 
Grebennikov, Sagdeev, Semchankau and Vasilevskii (2024).

We also obtain some lower bound estimates on the cardinality of the product set of two factorials modulo a prime. For instance, we prove that  if $N<p^{3/5},$ then
$$
\#\{m!n!\pmod p; \, 1\le m,n\le N\}\gg N^{1-o(1)}.
$$
The proof of this result is based on works of Banks and Shparlinski (2020), Cilleruelo and Garaev (2016), and the work of Katz and Shen (2008) related to the Ruzsa-Pl\"unnecke inequality from additive combinatorics.
\end{abstract}  

\subjclass[2020]{}

\maketitle


\section{Introduction and our results}

Throughout the paper, $p$ is a large prime number. 

There has been a much work in the literature related to the problem of distribution of factorials modulo $p$~\cite{CVZ, GaGa, GaGa1, GaHe, GaLuSh1, GaLuSh2,GaLuSh3,Garcia, Greb, KlMun, RSch}. 
It is a conjecture of Erd\"os (see,~\cite[{\bf
F11}]{RKG} and~\cite{RSch}) that about $p/e$ of the residue
classes modulo $p$ are missed by the sequence $n!$. If this were so,
the sequence $n!$ modulo $p$ should assume about $(1-1/e)p$ distinct
values, see~\cite{CVZ} for some results of probabilistic nature supporting  this conjecture. This in
turn would imply the representability of every residue class modulo
$p$ as a product of two factorials. 

The Wilson
theorem implies that
$$
\lambda!\cdot (p-\lambda)!\equiv \lambda \pmod p
$$
holds for any even $\lambda \in \{0, 2,\ldots, p-1\}$. The number of
such $\lambda$ is $(p+1)/2,$ and therefore by the pigeon-hole principle
every residue class modulo $p$ can be represented as
$$
m_1!n_1!+m_2!n_2!\pmod p.
$$
However, this argument does not apply to proving the existence of
representations involving factorials of integers of restricted size.
Some progress in this direction has been made in Garaev, Luca, Shparlinski~\cite{GaLuSh1,
GaLuSh2}. In particular, multiplicative character sums and double
exponential sums involving factorials have been estimated. These
estimates have been then applied to investigate various additive and
multiplicative congruences with factorials of integers in short
intervals.

In~\cite{GaLuSh2} it is
shown that any residue class $\lambda$ modulo $p$ can be represented
in the form
$$
\sum_{i=1}^{7}m_i!n_i! \equiv \lambda \pmod p
$$
with some positive integers $m_1, n_1, \ldots, m_7, n_7$ of the size
$O(p^{33/34}).$ This result was improved in~\cite{GaGa} in two directions, on the number of summands and the size of variables, that is, 
any residue class $\lambda$ modulo $p$ can be represented
in the form
$$
\sum_{i=1}^{5}m_i!n_i! \equiv \lambda \pmod p
$$
with some positive integers $m_1, n_1, \ldots, m_5, n_5$ of the size
$O(p^{27/28}).$ 

In the present paper we improve this result further to the following statement.
\begin{theorem}
\label{thm:main} Any residue class $\lambda$ modulo $p$ can be
represented in the form
$$
\sum_{i=1}^{5}m_i!n_i! \equiv \lambda \pmod p
$$
for some positive integers $m_1, n_1, \ldots, m_5, n_5$ with
$\max\limits_{1\le i\le 5}\{m_i, n_i\}\le p^{97/113+o(1)}.$ 
\end{theorem}

A challenging open question is to prove the existence of a fixed 
positive integer constant $\ell$ such that for every integer $\lambda$
the congruence
$$
n_1!+ \ldots + n_\ell!\equiv \lambda\pmod{p},
$$
has a solution in positive integers $n_1, \ldots, n_\ell$. An
estimate of the type
$$
\max_{\gcd(a, p) = 1} \left|\sum_{x=1}^{p}\exp(2\pi ia x!/p)\right|
\ll p^{1-c},
$$
for some constant $c>0,$ would completely solve this problem. The above mentioned conjecture of Erd\"os conditionally
implies that one can take $\ell=2.$

As an approximation to this open problem, it was shown in~\cite{GaLuSh2} that 
any residue class $\lambda$ modulo $p$ is representable in the form 
$$
m!n! +\sum_{i=1}^{47}n_i!\equiv \lambda\pmod p,
$$
with $\max\{m, n, n_1,\ldots, n_{47}\} \ll p^{1350/1351}.$
Our next result improves the exponent to $1300/1301.$

\begin{theorem}
\label{thm:main47} Any residue class $\lambda$ modulo $p$ can be
represented in the form
$$
m!n! +\sum_{i=1}^{47}n_i! \equiv \lambda \pmod p
$$
for some positive integers $m, n, n_1,\ldots, n_{47}$ with
$$
\max\{m, n, n_1,\ldots, n_{47}\} \ll p^{1300/1301}.
$$ 
\end{theorem}
 
It is interesting to note that there exists an integer $c_p$ such that the congruence
$$
x!+y!+c_pz!+c_pt!\equiv\lambda\pmod p, \quad 1\le x,y,z,t\le p,
$$
is solvable for any residue class $\lambda$ modulo $p,$ see~\cite{GaGa}. However,  the exact value of $c_p$ is not known.

\bigskip

Let 
$$
\cA_N=\{n!\pmod p, \quad 1\le n\le N \}.
$$
We consider the cardinality of the product set of two factorials
$$
\cA_N\cA_N= \{m!n!\pmod p,\, 1\le m,n\le N\}.
$$
The problem of obtaining lower bound estimates for the cardinality of the set $\cA_N\cA_N$ has been considered in several  papers, see~\cite{Garcia, GaGa1, GaLuSh1, Greb}. Such estimates are interesting by its own way, and also are useful in the study of multiplicative congruences with factorials. 

From the mentioned work of Grebennikov et. al.~\cite{Greb} it is known that 
$$
|\cA_N\cA_N|\gg p^{3/4},  \quad \text{if} \quad N>p^{7/8}\log p.
$$ 

It is also known that 
$$
|\cA_N\cA_N|\gg N, \quad \text{if} \quad N<p^{1/2},
$$
see~\cite[Section 5]{GaHe}.

Here, we shall prove the following lower bounds on the remaining range of the parameter $N,$ which improves some results
from~\cite{GaGa}.

\begin{theorem}
\label{thm:cardAA}
The following bound holds: 
\begin{equation*}
|\cA_N  \cA_N| \gg \begin{cases} 
N/(p^{1/8}\log p),         & \text{if}\    p^{29/40}\log p   \le N \le   p^{7/8}\log p,\\
\min\{p^{3/5}, N^{2/3}p^{1/6-o(1)}\},               & \text{if}\    p^{1/2}   \le N \le  p^{29/40}\log p.
\end{cases} 
\end{equation*}
\end{theorem}

In the range $N<p^{13/20}$  we can get a better bound.

\begin{theorem}
\label{thm:cardAAcg}
If $N<p^{3/5},$ then 
\begin{equation*}
|\cA_N  \cA_N| \gg N^{1-o(1)}.
\end{equation*}
\end{theorem}

The proof of Theorem~\ref{thm:cardAAcg} uses the following result on cardinality of product set, which is of independent interest.

\begin{theorem}
\label{thm:I times M}
Let $\mathcal{M}$ be arbitrary subset of $\mathbb{F}_p^*$ with $|\cM|=M$ elements, and let 
$\mathcal I$ be the interval
$$
\mathcal{I}=\{1,2,\ldots, N\}\pmod p\subset \mathbb{F}_p^*
$$
with $|\mathcal I|=N<p$ elements.
Then for the cardinality of the product set $\mathcal I\cdot\mathcal M$ we have the bound 
$$
|\mathcal I\cdot \mathcal{M}|\gg \min\Bigl\{\frac{p}{\log^2N},\, \frac{N^2}{\log ^2N},\, \frac{NM}{\log N}\Bigr\}. 
$$
\end{theorem}

\subsection{Notation and conventions}

Throughout the paper, we  use the notation
$$
\cA_N=\{n!\pmod p: 1\le n\le N\}.
$$
For simplicity, we will oftenly use $\cA$ instead of $\cA_N.$

We recall
that the notations $U=O(V),$ $U \ll V$ and  $V \gg U$  are
all equivalent to the statement that $|U| \le c V$ holds
with some constant $c> 0$. 

Any implied constants in symbols $O$, $\ll$
and $\gg$ may occasionally, where obvious, depend on the real parameter $\varepsilon >0,$ 
and on some other fixed constants.

We use $U< V p^{o(1)}$ to mean that for any epsilon $\varepsilon>0$, there exists $c=c(\varepsilon)>0$
such that $U < c V p^{\varepsilon}.$ Similarly, we use $U > Vp^{-o(1)}$ to mean that for 
any epsilon $\varepsilon>0$, there exists $c=c(\varepsilon)>0$ such that $U > c V p^{-\varepsilon}.$

Given a set $\cX$, we use $|\cX|$ to denote its cardinality.

We also use the abbreviation
$$
\ep(z) = \exp(2 \pi i z/p).
$$

\section{Lemmas}

Our proofs are based on recent results on distribution of factorials modulo $p$  obtained by Grebennikov,  Sagdeev, Semchankau and Vasilevskii~\cite{Greb}. 

\begin{lemma}
\label{lem:Greb}
Let 
$$
N\le p,\quad K =\frac{p}{N},\quad Q=\frac{N}{p^{1/2}(\log p)^2}.
$$ Then 
$$
|\cA_N/\cA_N|\gg \begin{cases} p, & \text{if}\ N\ge p^{7/8}\log p,\\
NQ^{1/3}(\log Q)^{-2/3}, & \text{if}\ p^{7/8}\log p\ge N\ge p^{4/5}(\log p)^{8/5},\\
NK^{1/2}, & \text{if}\ p^{4/5}(\log p)^{8/5}\ge N\ge p^{4/5}(\log p)^{4/5},
\\
NQ^{1/3}, & \text{if}\ p^{4/5}(\log p)^{4/5}\ge N\ge cp^{1/2}(\log p)^2,
\end{cases} 
$$
for some suitable positive absolute constant $c.$
\end{lemma}

Lemma~\ref{lem:Greb} is particular case of a more general result of~\cite[Theorem 1.3]{Greb}.

We  need the following statement from~\cite[Theorem 1]{GaLuSh2}.

\begin{lemma}
\label{lem2:GLSh2} Let $L$ and $N$ be integers with $0<L < L + N < p.$ 
Then for any fixed positive integer $\ell,$ the number $J_{\ell}(L,N)$ of 
solutions of the congruence
\begin{equation}
\label{eqn:sum n! = sum n!} 
\sum_{i=1}^{\ell}n_i!\equiv \sum_{i=\ell+1}^{2\ell}n_i!\pmod p,  \quad L+1\le n_1,\ldots, n_{2\ell}\le L+N,
\end{equation}
satisfies
$$
J_{\ell}(L,N)\ll N^{2\ell-1+1/(\ell+1)}.
$$
\end{lemma}

From Lemma~\ref{lem2:GLSh2} we have the following consequence.  
\begin{lemma}
\label{lem1:GLSh2}
Let $L$ and $N$ be integers with $0<L < L + N < p.$ Then the following estimate holds:
\begin{equation}
\label{eqn:doublesum} \max_{\gcd(a, p) = 1} \left|
\sum_{n=L+1}^{L+N} \sum_{x\in\cA} \ep(n! x)\right| \ll  |\cA|^{5/6}N^{7/8} p^{1/6}.
\end{equation}
\end{lemma}

Indeed, the Holder's inequality applied to the sum over $x$, gives 
\begin{eqnarray*}
\Bigl|\sum_{n=L+1}^{L+ N}\sum_{x\in\cA}\ep(an!x)\Bigr|^6\le |\cA|^5\sum_{x\in\cA}\Bigl|\sum_{n=L+1}^{L+ N}\ep(an!x)\Bigr|^6 \\ \le  |\cA|^5\sum_{x=0}^{p-1}\Bigl|\sum_{n=L+1}^{L+ N}\ep(an!x)\Bigr|^6 = p|\cA|^5 J_3(L, N),
\end{eqnarray*}
where $J_3(L, N)$ is the number of solutions of the congruence~\eqref{eqn:sum n! = sum n!} with $\ell =3.$ From Lemma~\ref{lem2:GLSh2} we have
$$
J_3(L, N)\ll N^{21/4},
$$
and the claim follows.

\bigskip

Let $\mathcal{M}$ be an arbitrary subset of $\mathbb{F}_p^*$ with $|\cM|=M$ elements, and $J(N,\mathcal{M})$ be the number of solution to the congruence 
$$ 
n_1 m_1 \equiv n_2 m_2  \pmod p, \quad 1 \leq n_1, n_2 \leq N, \quad m_1,m_2 \in \mathcal{M}.  
$$

The following result is due to Banks and Shparlinski~\cite[Theorem 2.1]{BS}. It will be used in the proof of Theorem~\ref{thm:cardAA}.

\begin{lemma}
\label{lem:BanksSh}
The following bound holds
\begin{equation*}
J(N,\mathcal{M}) \ll  \begin{cases} N^2 M^2 p^{-1} + N M p^{o(1)} , & \text{if}\ N\ge p^{2/3} ,\\
N^2 M^2 p^{-1} + N M^{7/4} p^{-1/4+o(1)} + M^2, & \text{if}\ N< p^{2/3} \; \text{and} \; M \ge p^{1/3},\\
N M p^{o(1)}  + M^2, & \text{if}\ N< p^{2/3} \; \text{and} \; M < p^{1/3}.
\end{cases} 
\end{equation*}
\end{lemma}
Using the well-known relationship between cardinality of a product set with the number of solutions of the corresponding multiplicative congruence, Lemma~\ref{lem:BanksSh} implies the following consequence.
\begin{cor}
\label{cor:BanksSh} Let $\cI$ be the interval
$$
\mathcal{I}=\{1,2,\ldots, N\}\pmod p\subset \mathbb{F}_p^*
$$
with $|\cI|=N<p$ elements.  If $N\ge p^{2/3},$ then
$$
|\cI\cdot \cM|\gg \min\{p, N M p^{-o(1)}\}.
$$
If $p^{1/2}\le N\le p^{2/3},$ then
$$
|\cI\cdot \cM| \gg \min\{p, N M^{1/4} p^{1/4-o(1)}, NMp^{-o(1)}\}.
$$
\end{cor}

In the proof of Theorem~\ref{thm:cardAA} and Theorem~\ref{thm:cardAAcg} we will also use the following two statements from additive combinatorics. The first one is the  multiplicative form of a result of Katz and Shen~\cite[Corollary 1.5]{KS} related to the Ruzsa-Pl\"unnecke inequality. The second is a particular case of Ruzsa's triangle inequality, see~\cite{Nat, R1, R2, TV}. All sets are considered to be non-empty.
\begin{lemma}
\label{lem:KatzShen}
Let $X,B_1,\ldots , B_k$ be any subsets of $\mathbb{F}^*_p$. Then there exists a subset $X' \subset X$ with $|X'| > 0.5 |X|$ so that
\begin{equation}
\label{eqn:K-S}
| X' B_1 \cdots B_k | \ll \frac{|X B_1| \cdots |X B_k| }{|X|^{k-1}}.
\end{equation}
\end{lemma}

\begin{lemma}
\label{lem:Ruzsa}
For any  subsets $X, Y, Z\subset \mathbb{F}^*_p$ we have
$$
|X/Y| \le \frac{|XZ|\cdot|ZY|}{|Z|}.
$$
\end{lemma}

In the proof of Theorem~\ref{thm:I times M} we will  use the following result of Cilleruelo and Garaev~\cite{CiGa}.

\begin{lemma}
\label{lem:CG}
For any $s_0\in \F_p$, the number of
solutions of the congruence
\begin{equation}
\label{eqn:x=sy}
x\equiv s_0y\pmod p;\quad x,y\in \mathbb{N},\quad x\le X,\quad y\le Y,\quad \gcd(x,y)=1,
\end{equation}
is bounded by $O(1+XYp^{-1}).$
\end{lemma}

\section{Proof of Theorem~\ref{thm:main}}

We can assume that $p^{4/5}(\log p)^{8/5}<p^{97/113}\log p<N<p^{7/8}.$ In particular, from Lemma~\ref{lem:Greb}, it follows that
\begin{equation}
\label{eqn: A/A gg N43p16logp}
|\cA/\cA|\gg \frac{N^{4/3}}{p^{1/6}(\log p)^{4/3}}.
\end{equation}

Let 
$$
\{(x_j,y_j): \, 1\le j\le |\cA/\cA|\} \subset \cA\times \cA
$$ 
be the set of all $|\cA/\cA|$ pairs $(x_j,y_j)$ for which  all the elements $x_j/y_j\pmod p$
are pairwise distinct.

Let $J$ be the number of solutions of the congruence 
\begin{equation}
\label{eqn: firstmain}
x_j(a_1+a_2)+y_j(a_3+a_4) +n!a_5\equiv \lambda \pmod p,
\end{equation}
in variables $j,a_1,a_2,a_3,a_4, a_5, n$ satisfying 
$$
j\le |\cA/\cA|, \quad a_1,a_2,a_3,a_4, a_5\in \cA, \quad n\le N. 
$$
It suffices to prove that $J>0.$ 

We express $J$ in terms of exponential sums. 
$$
J=\frac{1}{p}\sum_{a=0}^{p-1}\sum_{n, a_5}\ep(an!a_5)\sum_{j}\sum_{a_1,a_2,a_3,a_4}
\ep(a(x_j(a_1+a_2) +y_j(a_3+a_4)-\lambda)).
$$
Picking up the term that corresponds to $a=0$ and then applying~\eqref{eqn:doublesum}
to the sum over $n,a_5,$ we get
\begin{equation}
\label{eqn: firstJ}
J= \frac{|\cA/\cA|\cdot|\cA|^5 N}{p} +R,
\end{equation}
where
\begin{eqnarray*}
|R|\ll \frac{|\cA|^{5/6}N^{7/8} p^{1/6}}{p}\sum_{a=1}^{p-1}\sum_{j}\Bigl|\sum_{a_1,a_2,a_3,a_4}
\ep(a(x_j(a_1+a_2) +y_j(a_3+a_4))\Bigr|.
\end{eqnarray*}
Since
\begin{eqnarray*}
&&\Bigl|\sum_{a_1,a_2,a_3,a_4\in\cA}
\ep(a(x_j(a_1+a_2) +y_j(a_3+a_4))\Bigr|=\\
&&\qquad \Bigl|\sum_{a_1,a_3\in\cA}\ep(a(x_ja_1+y_ja_3)\Bigr|\Bigl|\sum_{a_2,a_4\cA}
\ep(a(x_ja_2 +y_ja_4)\Bigr|=\\
&&\qquad\qquad\Bigl|\sum_{a_1,a_3\in \cA}
\ep(a(x_ja_1+y_ja_3))\Bigr|^2,
\end{eqnarray*}
we get that
$$
|R| \ll
\frac{|\cA|^{5/6}N^{7/8} p^{1/6}}{p}\sum_{a=0}^{p-1}\sum_{j\le |\cA/\cA|}\Bigl|\sum_{a_1,a_3\in \cA}
\ep(a(x_ja_1+y_ja_3))\Bigr|^2.
$$
Note that the quantity
$$
\frac{1}{p}\sum_{a=0}^{p-1}\sum_{j\le |\cA/\cA|}\Bigl|\sum_{a_1,a_3\in \cA}
\ep(a(x_ja_1+y_ja_3))\Bigr|^2=J_0,
$$
where $J_0$ is the number of solutions of the congruence
$$
x_ja_1+y_ja_3 \equiv x_ja_2+y_ja_4\pmod p, \, j\le |\cA/\cA|, a_1,a_2,a_3,a_4\in \cA.
$$
Hence, 
\begin{eqnarray*}
|R|\ll  |\cA|^{5/6}N^{7/8} p^{1/6} J_0.
\end{eqnarray*}
We rewrite the equation in the form
$$
\frac{x_j}{y_j}(a_1-a_2)\equiv a_4-a_3\pmod p, \, j\le |\cA/\cA|, a_1,a_2,a_3,a_4\in \cA.
$$

If $a_4\equiv a_3\pmod p$ then $a_1\equiv a_2\pmod p$ and $j\le |\cA/\cA|$ can be arbitrary. This case contributes to $J_0$ a quantity of order $|\cA/\cA||\cA|^2.$

If $a_4\not \equiv a_3\pmod p,$ then $a_1\not \equiv a_2\pmod p.$ In this case, for fixed $a_1,a_2,a_3,a_4\in \cA$ we get $x_j/y_j$ uniquely determined, and thus we get $j$ determined (recall that all the $|\cA/\cA|$ elements $x_j/y_j$ are distinct modulo $p$). Hence, this case contributes to $J_0$ a quantity not greater than $|\cA|^4.$

Therefore, since $|\cA/\cA|\le |\cA|^2,$  
$$
J_0\le |\cA/\cA|\cdot |\cA|^2 + |\cA|^4 \le 2|\cA|^4.
$$
Inserting this to the above estimate of $R,$ gives
$$
|R|\ll |\cA|^{5/6}N^{7/8} p^{1/6}|\cA|^4.
$$
Thus, from~\eqref{eqn: firstmain} we get that
$$
J= \frac{|\cA/\cA|\cdot|\cA|^5 N}{p}\Bigl(1+O\Bigl(\frac{p^{7/6}}{|\cA/\cA|\cdot|\cA|^{1/6}N^{1/8}}\Bigr)\Bigr).
$$

The result now follows from~\eqref{eqn: A/A gg N43p16logp} and $|\cA|\ge |\cA/\cA|^{1/2}$.

\section{Proof of Theorem~\ref{thm:main47}}

Since $N>p^{7/8}\log p,$ from Lemma~\ref{lem:Greb} it follows that 
$$
|\cA|\ge |\cA/\cA|^{1/2}\gg p^{1/2}.
$$

We shall prove that for some suitable constant $c>0,$ if $N>cN^{1300/1301}$ then the congruence
$$
n! x+\sum_{i=1}^{47}n_i!\equiv \lambda \pmod p,
$$
has a solution in positive integers with 
$$
\max\{n, n_1,\ldots, n_{47}\}\le N,\quad x\in \cA.
$$ 
Expressing the number $J$ of solutions of this congruence in terms of exponential sums,
we get
$$
J=\frac{1}{p}\sum_{a=0}^{p-1}\sum_{n\le N}\sum_{x\in\cA}\,
\sum_{n_1\le N,\ldots,n_{47}\le N}\ep(a(n!x+n_1!+\ldots+n_{47}!-\lambda)).
$$
Separating the term that corresponds to $a=0,$ and following the standard procedure, we get that
\begin{equation}
\label{eqn:J=N48Ap+R}
J=\frac{N^{48}|\cA|}{p} +R, 
\end{equation}
where
\begin{equation}
\label{eqn:R47}
|R|\le \frac{1}{p}\sum_{a=1}^{p-1}\Bigl(
\Bigl|\sum_{n\le N}\sum_{x\in\cA}\ep(an!x)\Bigr|\cdot \Bigl|\sum_{n\le N}\ep(an!)\Bigr|^{47}\Bigr).
\end{equation}
Taking  into account~\eqref{eqn:doublesum}, we get
$$
|R|\ll  \frac{p^{1/6} |\cA|^{5/6}N^{7/8}}{p}\, \sum_{a=0}^{p-1} \Bigl|\sum_{n\le N}\ep(an!)\Bigr|^{47}.
$$
As in~\cite{GaLuSh2}, the Cauchy-Schwarz inequality gives
\begin{eqnarray*}
&&\frac{1}{p}\sum_{a=0}^{p-1} \Bigl|\sum_{n\le N}\ep(an!)\Bigr|^{47}\le\\
&&\qquad \Bigl(\frac{1}{p}\sum_{a=0}^{p-1} \Bigl|\sum_{n\le N}\ep(an!)\Bigr|^{46}\Bigr)^{1/2}
\cdot 
\Bigl(\frac{1}{p}\sum_{a=0}^{p-1} \Bigl|\sum_{n\le N}\ep(an!)\Bigr|^{48}\Bigr)^{1/2}=\\
 &&\qquad \qquad \Bigl(J_{23}(0,N)\cdot J_{24}(0,N)\Bigr)^{1/2}.
\end{eqnarray*}

Using Lemma~\ref{lem2:GLSh2} with $\ell =23$ and $\ell=24$, we obtain that
$$
\frac{1}{p}\sum_{a=0}^{p-1} \Bigl|\sum_{n\le N}\ep(an!)\Bigr|^{47}\ll
N^{46+49/1200}.
$$
Hence, for some absolute constant $c_1>0,$ we have that
$$
|R|\le c_1p^{1/6}|\cA|^{5/6}N^{7/8}N^{46+49/1200}. 
$$
Since $|\cA|\ge c_2p^{1/2}$ for some absolute positive constant $c_2,$ it follows from~\eqref{eqn:J=N48Ap+R}, that
\begin{eqnarray*}
J\ge \frac{N^{48}|\cA|}{p}- c_1p^{1/6}|\cA|^{5/6}N^{7/8}N^{46+49/1200} \\
= \frac{N^{48}|\cA|}{p} \Bigl(1-\frac{c_1p^{7/6}}{|\cA|^{1/6}N^{1301/1200}}\Bigr)\\
>\frac{N^{48}|\cA|}{p} \Bigl(1-\frac{c_1p^{13/12}}{c_2^{1/6}N^{1301/1200}}\Bigl)>0,
\end{eqnarray*}
assuming that $N>cp^{1300/1301}$ with a sufficiently large constant $c.$

\section{Proof of Theorem~\ref{thm:I times M}}

Let $\mathcal{P}_N$ be the set of prime numbers not exceeding $N.$ Denote by $J$ the number of solutions of the congruence
\begin{equation*}
\label{eqn:qm=rn}
q m\equiv rn\pmod p, \quad q, r\in\cP_N,\quad m,n\in\cM.
\end{equation*}
Given $m,n\in\cM,$ denote by $J(m,n)$ the number of solutions of the congruence
\begin{equation}
\label{eqn:qm=rn with m and n fixed}
q m\equiv rn\pmod p, \quad q, r\in\cP_N.
\end{equation}
Then,
$$
J=\sum_{m\in\cM, n\in\cM} J(m,n).
$$
 If $m=n$ then $q=r$ and we have at most $\pi(N)M\ll NM/\log N$
contribution from this case. Thus,
\begin{equation}
\label{eqn:sum over m and n of J(m,n)}
J\ll \frac{NM}{\log N} +\sum_{\substack{m\in\cM, n\in\cM\\m\not=n}} J(m,n).
\end{equation}
From~\eqref{eqn:qm=rn with m and n fixed} we see that $J(m,n)$ is equal to the number of solutions of the congruence
$$
q \equiv (nm^*)r\pmod p, \quad q, r\in\cP_N,
$$
where $m^*$ is the multiplicative inverse of $m\pmod p.$ Since $m\not=n,$ we have $q\not=r$ and therefore, $\gcd(q,r)=1.$ It follows that for $m\not=n$ the number $J(m,n)$ 
does not exceed the number of solutions of the congruence
$$
x\equiv s_0y\pmod p,\quad x,y\in \mathbb{N},\quad x\le N,\quad y\le N,\quad \gcd(x,y)=1,
$$ 
where $s_0=nm^*.$ Hence, by Lemma~\ref{lem:CG}, for $m\not=n$ we have
$$
J(m,n)\ll 1+\frac{N^2}{p}.
$$
Thus, from~\eqref{eqn:sum over m and n of J(m,n)}, we get
$$
J\ll \frac{NM}{\log N} + M^2+\frac{N^2M^2}{p}.
$$

Using the relationship between the cardinality of a product set and the number of solutions of the corresponding multiplicative congruence, we conclude that
\begin{eqnarray*}
|\cI\cdot \cM|\ge |\cP_N\cdot \cM|\gg \frac{|\cP_N|^2\cdot |M|^2}{J} \gg 
\frac{(\pi(N))^2M^2}{\frac{NM}{\log N}+M^2+\frac{M^2N^2}{p}} \gg \\
\min\Bigl\{\frac{NM}{\log N}, \, \frac{N^2}{\log^2N},\, \frac{p}{\log^2N}\Bigr\}.
\end{eqnarray*}

\section{Proof of Theorem~\ref{thm:cardAA} and~Theorem~\ref{thm:cardAAcg}}

From Lemma~\ref{lem:Ruzsa} we have that
$$ 
|\cA / \cA| \leq \frac{|\cA \, \cA|^2}{| \cA |}.
$$
Since $|\cA|\ge |\cA/\cA|^{1/2}$, we get
\begin{equation}
\label{eqn:AA3/4}
|\cA \, \cA| \geq |\cA / \cA|^{3/4}.
\end{equation}
In particular, from Lemma~\ref{lem:Greb} we get the first inequality of Theorem~\ref{thm:cardAA}. 

Next, in view of Lemma~\ref{lem:KatzShen}, there exists a subset, $\cA' \subset \cA  $ such that $|\cA' |\geq 0.5 |\cA|$ and 
$$
| \cA' \, \cA /  \cA |  \leq \frac{ | \cA \, \cA  | \cdot | \cA / \cA    |}{ | \cA |} ,
$$
implying 
$$
|\cA \, \cA| \cdot | \cA / \cA |^{1/2}  \gg  |\cA' \, \cA /  \cA   | .
$$
Since $\cI:=\{1,2,\ldots, N\}\subset \cA/\cA$, we get that
$$
|\cA' \, \cA /  \cA |\ge |\cI\cdot \cA'|,
$$
whence
\begin{equation}
\label{eqn:KShtoFactorials}
|\cA \, \cA| \cdot | \cA / \cA |^{1/2}  \gg  |\cI\cdot \cA'|.
\end{equation}
We recall that 
\begin{equation}
\label{eqn:A'via A/A}
|\cA'|\ge 0.5|\cA|\ge 0.5|\cA/\cA|^{1/2}.
\end{equation}

To prove the second inequality of Theorem~\ref{thm:cardAA}, we will consider two cases.

\bigskip

{\bf Case 1.} Let $ p^{2/3}\le  N<p^{29/40}\log p.$
From Corollary~\ref{cor:BanksSh}, we have that
$$
|\cI\cdot \cA'|\gg \min \{p, N|\cA'|p^{-o(1)}\}\gg \min \{p, N|\cA/\cA|^{1/2}p^{-o(1)}\}.
$$

In the case $|\cI\cdot \cA'|\gg p$, from~\eqref{eqn:KShtoFactorials} we get
$$
|\cA \, \cA| \cdot | \cA / \cA |^{1/2}  \gg p.
$$
Taking this into account in~\eqref{eqn:AA3/4}, we obtain that
$$
|\cA\,\cA|^{5/2}\gg p^{3/2},
$$
which implies the desired bound $|\cA\,\cA|\gg p^{3/5}.$

In the case $|\cI\cdot \cA'|\gg  N|\cA/\cA|^{1/2}p^{-o(1)},$ from~\eqref{eqn:KShtoFactorials} we have that
$$
|\cA \, \cA| \cdot | \cA / \cA |^{1/2} \gg N|\cA/\cA|^{1/2}p^{-o(1)}.
$$
Therefore,
$$
|\cA \, \cA|\gg Np^{-o(1)}\gg p^{2/3-o(1)}\gg p^{3/5},
$$
and we are done in the Case 1.

\bigskip

{\bf Case 2.} Let $p^{1/2}\le  N<p^{2/3}.$ Our aim is to prove the bound
$$
|\cA\,\cA|\gg \min \{p^{3/5}, N^{2/3}p^{1/6-o(1)}\}.
$$

From Corollary~\ref{cor:BanksSh} and~\eqref{eqn:A'via A/A} we have 
$$
|\cI\cdot \cA'| \gg \min\{p, N |\cA/\cA|^{1/8} p^{1/4-o(1)}, N|\cA/\cA|^{1/2}p^{-o(1)}\}.
$$

If $|\cI\cdot \cA'|\gg p,$ then as before from~\eqref{eqn:AA3/4} and~\eqref{eqn:KShtoFactorials} we get 
$|\cA \, \cA|\gg p^{3/5}.$

If $|\cI\cdot \cA'| \gg  N |\cA/\cA|^{1/8} p^{1/4-o(1)}$, then from~\eqref{eqn:KShtoFactorials} we see that
$$
|\cA\,\cA|\cdot |\cA/\cA|^{3/8}\gg Np^{1/4-o(1)}
$$
which together with~\eqref{eqn:AA3/4} gives
$$
|\cA\,\cA|\gg N^{2/3}p^{1/6-o(1)}.
$$

If
$$
|\cI\cdot \cA'| \gg  N|\cA/\cA|^{1/2}p^{-o(1)},
$$
then from~\eqref{eqn:KShtoFactorials} we get
$$
|\cA\,\cA|\gg Np^{-o(1)}\gg N^{2/3}p^{1/6-o(1)}
$$
finishing the proof of Theorem~\ref{thm:cardAA}.

In order to prove Theorem~\ref{thm:cardAAcg}, we can assume that $N>p^{1/2}$ (recall that the statement is known for $N<p^{1/2}$). Hence, from~Theorem~\ref{thm:I times M} and~\eqref{eqn:A'via A/A} we get that
$$
|\cI\cdot \cA'|\gg \min\Bigl\{\frac{p}{\log^2N},\,  \frac{N|\cA/\cA|^{1/2}}{\log N}\Bigr\}. 
$$
In view of~\eqref{eqn:KShtoFactorials} we have 
$$
|\cA \, \cA| \cdot | \cA / \cA |^{1/2}  \gg  \min\Bigl\{\frac{p}{\log^2N},\,  \frac{N|\cA/\cA|^{1/2}}{\log N}\Bigr\}. 
$$

If 
$$
|\cA \, \cA| \cdot | \cA / \cA |^{1/2}  \gg  \frac{p}{\log^2N},
$$
then together with~\eqref{eqn:AA3/4}, we get that
$$
|\cA\, \cA|\gg \frac{p^{3/5}}{(\log N)^{6/5}}\gg N^{1-o(1)}. 
$$

If
$$
|\cA \, \cA| \cdot | \cA / \cA |^{1/2}  \gg  \frac{N|\cA/\cA|^{1/2}}{\log N}, 
$$
we get 
$$
|\cA \, \cA|  \gg  \frac{N}{\log N}>N^{1-o(1)},
$$
which finishes the proof of Theorem~\ref{thm:cardAAcg}.

\end{document}